\documentclass[12pt]{article}

\usepackage{amsmath,amssymb,amsthm}
\usepackage[a4paper,scale=0.8]{geometry}
\usepackage{graphicx,color}

\theoremstyle{remark}

\newcommand{\bs}{\boldsymbol}
\def\xb{{\mathbf x}}
\def\Pad{\operatorname{Pad}}

 \def\xb{{\mathbf x}}

 \def\CO{{\mathcal O}}

\title{Polynomial approximation and cubature \\ at approximate Fekete and Leja points of the cylinder}

\author{
Stefano De Marchi\thanks{University of Padua (Italy), {\tt
demarchi@math.unipd.it}},
Martina Marchioro\thanks{University of Padua (Italy), {\tt
martina.marchioro@gmail.com}} and
Alvise Sommariva\thanks{University of Padua (Italy), {\tt
alvise@math.unipd.it}}
}

\date{\today}

\begin{document}

\maketitle

\begin{abstract}
The paper deals with polynomial interpolation, least-square
approximation and cubature of functions defined on the rectangular
cylinder, $K=D\times [-1,1]$, with $D$ the unit disk. The nodes used
for these processes are the {\it Approximate Fekete Points} (AFP)
and the {\it Discrete Leja Points} (DLP) extracted from suitable
{\it Weakly Admissible Meshes} (WAMs) of the cylinder. 
From the
analysis of the growth of the Lebesgue constants, approximation and
cubature errors, we show that the AFP and the DLP extracted from WAM
are good points for polynomial approximation and numerical
integration of functions defined on the cylinder.
\end{abstract}

\section{Introduction}
Locating good points for multivariate polynomial
approximation, in particular interpolation, is an open challenging
problem, even in standard domains.

One set of points that is always
good, in theory, is the so-called {\em Fekete points}. They are
defined to be those points that maximize the (absolute value of the)
Vandermonde determinant on the given compact set. However, these are
known analytically only in a few instances (the interval and the
complex circle for univariate interpolation, the cube for tensor
product interpolation), and are very difficult to compute, requiring
an expensive and numerically challenging nonlinear multivariate
optimization.

{\em Admissible Meshes} (shortly AM), introduced by Calvi and
Levenberg in \cite{CL08}, are sets of points in a given compact
domain $K \subset \mathbb{R}^d$ which are nearly optimal for
least-squares approximation, and contain interpolation points that
distribute asymptotically as Fekete points of the domain. This
theory has given new insight to the (partial) solution of the
problem of extracting good interpolation point in dimension $d>1$.
In all practical applications, instead of AM, people look for
low-cardinality admissible meshes, called {\em Weakly Admissible
Meshes} or WAM (cf. the recent survey \cite{BDMSV09_1}).

The extremal sets of our interest are the {\em Approximate Fekete Points} (AFP) and {\it Discrete
Leja Points} (DLP). As described in \cite{BDMSV09,SV09-1}, AFP and DLP
can be easily computed by using basic tools of numerical linear algebra.
In practice, (Weakly) Admissible Meshes and Discrete Extremal Sets allow us to replace a
continuous compact set by a discrete version, that is ``just as
good'' for all practical purposes.

In this paper, we focus on AFP and DLP extracted from Weakly
Admissible Meshes of the rectangular cylinder $K=D \times [-1,1]$,
with $D$ the unit disk. These points are then used for computing
interpolants, least-square approximants and cubatures of functions
defined on the cylinder. We essentially provide AFP and DLP of the
cylinder that, to our knowledge, have never been investigated so
far.

\section{Weakly Admissible Meshes: definitions, properties and construction}

Given a {\em polynomial determining} compact set $K\subset
\mathbb{R}^d$ or $K\subset \mathbb{C}^d$ (i.e., polynomials
vanishing there are identically zero), a Weakly Admissible Mesh
(WAM) is defined in \cite{CL08} to be a sequence of discrete subsets
${A}_n\subset K$ such that
\begin{equation} \label{WAM}
\|p\|_K\leq C(A_n) \|p\|_{{A}_n}\;,\;\;\forall
p\in \mathbb{P}_n^d(K)
\end{equation}
$\mathbb{P}_n^d(K)$ being the set of $d$-variate polynomials of degree at most $n$ on $K$,
where both ${\rm card}({A}_n)\geq N:=\mbox{dim}(\mathbb{P}_n^d(K))$ and $C(A_n)$
grow at most {\em polynomially} with $n$, i.e. {\tt card}$(A_n) \le c \;n^s$, for some fixed
$s\in \mathbb{N}$ depending only on $K$. When $C({A}_n)$ is
bounded we speak of an Admissible Mesh (AM). We use
the notation $\|f\|_K=\sup_{x\in K}|f(x)|$ for $f$ a bounded
function on the compact $K$.

\noindent WAMs enjoy the following ten properties (already
enumerated in \cite{BCLSV09} and proved in \cite{CL08}): \vskip0.3cm
\noindent
{\bf P1:} $C({A}_n)$ is invariant under affine mapping\\
{\bf P2:} any sequence of unisolvent interpolation sets whose
Lebesgue constant grows at most polynomially with $n$ is a WAM,
$C({A}_n)$ being the Lebesgue constant itself\\
{\bf P3:} any sequence of supersets of a WAM whose cardinalities
grow polynomially with $n$ is a WAM with the same constant
$C({A}_n)$\\
{\bf P4:} a finite union of WAMs is a WAM for the corresponding
union of compacts, $C({A}_n)$ being the maximum of the
corresponding constants\\
{\bf P5:} a finite cartesian product of WAMs is a WAM for the
corresponding product of compacts, $C({A}_n)$ being the
product
of the corresponding constants\\
{\bf P6:} in $\mathbb{C}^d$ a WAM of the boundary $\partial K$ is
a WAM of $K$ (by the maximum principle)\\
{\bf P7:} given a polynomial mapping $\pi_s$ of degree $s$, then
$\pi_s({A}_{ns})$ is a WAM for $\pi_s(K)$ with constants
$C({A}_{ns})$ (cf. \cite[Prop.2]{BCLSV09})\\
{\bf P8:} any $K$ satisfying a Markov polynomial inequality like
$\|\nabla p\|_K\leq Mn^r \|p\|_K$ has an AM with
${O}(n^{rd})$
points (cf. \cite[Thm.5]{CL08})\\
{\bf P9:} least-squares polynomial approximation of $f\in C(K)$: the
least-squares polynomial ${L}_{{A}_n} f$ on a WAM is
such that
$$
\|f-{L}_{{A}_n} f\|_K \lessapprox
C({A}_n)\sqrt{{\rm
card}({A}_n)}\,\min{\{\|f-p\|_K,\,p\in \mathbb{P}_n^d(K)\}}
$$
(cf. \cite[Thm.1]{CL08})\\
{\bf P10:} Fekete points: the Lebesgue constant of Fekete points
extracted from a WAM can be bounded like $\Lambda_n\leq NC({A}_n)$
(that is the elementary classical bound of the continuum Fekete
points times a factor $C({A}_n)$); moreover, their asymptotic
distribution is the same of the continuum Fekete points, in the
sense that the corresponding discrete probability measures converge
weak-$\ast$ to the pluripotential equilibrium measure of $K$ (cf.
\cite[Thm.1]{BCLSV09}). Pluripotential theory has been widely
studied by M. Klimek in the monograph \cite{KL92}, to which
interested readers should refer for more details. \vskip0.3cm
\noindent It is worth noticing that in the very recent papers
\cite{K10, PV10}, the authors have provided new techniques for
finding admissible meshes with low cardinality, by means of
analytical transformations of domains. \vskip 0.2in \noindent
Examples of WAMs can be found in \cite{BDMSV09, BDMSV09_1}. Here, we
simply recall some one dimensional and two dimensional WAMs.
\begin{enumerate}
\item The set
$$C_n=\{\cos(k\pi/n), \;\; k=0,\ldots,n\}$$
of $n+1$ {\it Chebyshev-Lobatto} points for the interval $I=[-1,1]$, is a
one-dimensional WAM of degree $n$ with $C(A_n)=\CO(\log \,n)$ and {\tt card}$(C_n)=n+1$.
This follows from property {\bf P2}.
\item The set $\Pad_n,\; n\ge0$ of the {\em Padua points} of degree $n$
of the square $Q=[-1,1]^2$ is the set defined as follows (cf.
\cite{BCDMVX06})
\begin{equation} \label{padua1}
  \Pad_n=\{ \xb_{k,j} =  (\xi_k,\eta_j), \quad 0 \leq k\leq n, \quad
         1\leq j\leq \lfloor \tfrac{n}{2} \rfloor +1 \},
\end{equation}
where
\begin{equation} \label{padua2}
\xi_k=\cos{\frac{k \pi}{n}}, \qquad \eta_j=\left\{ \begin{array}
{lr}
   \cos{\frac{2j-1}{n+1}\pi}, & \;k \;\mbox{even}\\  &  \\
\cos{\frac{2j-2}{n+1}\pi},  & \;k \;\mbox{odd}
\end{array}
\right.
\end{equation}
Notice that here we refer to the first family of Padua points.
$\Pad_n$ is then a two-dimensional WAM with $C(\Pad_n)=\CO(\log^2
\,n)$ and {\tt card}$(\Pad_n)=(n+1)(n+2)/2$. This is a consequence
of property {\bf P2} since, as shown in \cite{BCDMVX06}, the Padua
points are a unisolvent set for polynomial interpolation in the
square with minimal order of growth of their Lebesgue constant, i.e.
$\CO(\log^2 \,n)$.
\item The sequence of {\em polar symmetric grids} $A_n=\{ (r_i\cos \theta_j, r_i \sin \theta_j)\}$ with
the radii and angles defined as follows
\begin{equation} \label{eq1}
(r_i,\theta_j)_{i,j}=\left\{ \cos(i\pi/n), \; 0\le i \le n \right\}
\times \left\{
 {j\pi \over n+1}, 0 \le j\le n \right \}
\end{equation}
are WAMs for the closed unit disk $D=\{{\bf x} \,: \, \|{\bf
x}\|_2\le 1\}$, with constant $C(A_n)=\CO(\log^2 n)$ and cardinality
{\tt card}$(A_n)=(n+1)^2$ for odd $n$ and {\tt card}$(A_n)=n^2+n+1$
for even $n$ (cf. \cite[Prop. 1]{BSV09}). Moreover, since these WAMs
contain the Chebyshev-Lobatto points of the vertical diameter
$\theta = \pi/2$ only for $n$ odd (whereas it always contains the
Chebyshev-Lobatto points of the horizontal diameter $\theta = 0$),
and thus is not invariant under rotations by an angle $\pi/2$. Hence
in order to have a WAMs on the disk invariant by rotations of
$\pi/2$, we have to modify the choice of radii and angles in
(\ref{eq1}) as follows
\begin{equation}
(r_i,\theta_j)_{i,j}=\left\{ \cos(i\pi/n), \; 0\le i \le n \right\}
\times \left\{
 {j\pi \over n+2}, 0 \le j\le n+1 \right \}, n \; even
\end{equation}
In this way the obtained WAM is now invariant with {\tt card}$(A_n)=(n+1)^2$ also for
$n$ even.
\end{enumerate}
\subsection{Three dimensional WAMs of the cylinder}
We restrict ourselves to the rectangular cylinder with unitary radius and height the interval
[-1,1], that is $K=D \times [-1,1]$, where as above, $D$ is the closed unit disk.

We considered two meshes: the first one uses a symmetric polar grid
in the disk $D$ and Chebyshev-Lobatto points along $[-1,1]$; the
second one uses Padua points on the $(x,z)$ plane and equispaced
points along the circumference of $D$.

\subsubsection{The first mesh: WAM1}
We consider the set
$$
A_n=\{ (r_i \cos\theta_j, r_i \sin\theta_j,z_k)\}
$$
with $-1\leq r_i\leq 1, \; 0\le i \le n$ and  $0\leq \theta_j \leq \pi, \; 0\le j\le n$ that is
\begin{eqnarray*}
\{(r_i, \theta_j, z_k)\}_{i,j,k} & = & \left\{ \cos\left( \frac{i
\pi}{n}\right), 0\leq i \leq n \right\}\times\left\{ \left.
\begin{array} {ll}
\frac{j \pi}{n+2}, 0\leq j \leq n+1,\,\; n\, \mbox{even}\\
\frac{j \pi}{n+1}, 0\leq j \leq n,\,\; n \,\mbox{odd}
\end{array}
\right.
\right\}\\
& \times & \left\{\cos\left( \frac{k \pi}{n}\right), 0\leq k \leq n  \right\}
\end{eqnarray*}
The cardinality of $A_n$, both for $n$ even and $n$ odd, is $(n+1)^3$.
Indeed, let us consider first the case of $n$ even.
The points on the disk, subtracting 
the repetitions of the center, which are $n+2-1$, are
$(n+1)^2$. All these points are then multiplied by
the corresponding $n+1$ Chebyshev-Lobatto points
along the third axis $z$, giving the claimed cardinality.

\noindent When $n$ is odd, there are no coincident points, thus we
have $(2n+2)(n+1)/2=(n+1)^2$ points on the disk. Then, considering
the $n+1$ Chebyshev-Lobatto points along the third axis, we get the
claimed results.

\noindent Finally, the set $A_n$ so defined, is a WAM since it is
the cartesian product of a two dimensional WAM (the points on the
disk) and the one dimensional WAM of the Chebyshev-Lobatto points.
The property {\bf P5} gives the constant $C(A_n)=\CO(\log^3 \,n)$
(see Figure \ref{Fig1} for the case $n=5$).
\subsubsection{The second mesh: WAM2}
This discretization is obtained by taking the Padua points $\Pad_n$
on the plane $(r,z)$, rotated $n$ times along $z$-axis by a constant
angle $\theta=\pi/(n+1)$. In this way, along the bottom
circumference of the cylinder, we obtain $2n+2$ equispaced points.
This is due to the fact that the points with coordinates $(-1,0)$
and $(1,0)$ are Padua points. In details, the mesh is the set
$$
A_n=\{ (r_i \cos\theta_j, r_i \sin\theta_j,z_k)\}
$$
with $-1\leq r_i\leq 1,\; 0 \le i \le n$ and  $0\leq \theta_j \leq \pi,\; 0\le j \le n$, that is
\begin{eqnarray*}
\{(r_i, \theta_j, z_k)\}_{i,j,k} & = & \left\{ \cos\left( \frac{i \pi}{n}\right), 0\leq i \leq n \right\}\times\left\{\frac{j \pi}{n}, 0\leq j \leq n+1\right\}\\
& \times & \left\{\cos\left( \frac{k \pi}{n+1}\right), 0\leq k \leq n+1 \,\,
\left.
\begin{array} {rl}
k & \mbox{odd when $i$ is even}\\
k & \mbox{even when $i$ is odd}
\end{array}
\right. \right\}.
\end{eqnarray*}
This mesh has cardinality $\CO\left({n^3 \over 2}\right)$. In fact, when $n$ is even, the points are ${(n+1)(n+2) \over 2} (n+1)$ from which
we have to subtract the repetitions $({n \over 2}+1)\,n$, corresponding to the ${n\over 2} +1$ Padua points with abscissa $x=0$ counted $n$ times.
Then, the points so generated are $(n^2+n+1){(n+2) \over 2}$. On the contrary, when $n$ is
odd, there are no intersections and so the total number of points is ${(n+1)^2 (n+2) \over 2}$ (see Figure \ref{Fig1} for the case $n=5,6$).

We now prove that this mesh is indeed a WAM.
To this aim, consider a generic polynomial of degree at most $n$ defined on the cylinder
$p(x,y,z)=p(r \cos \theta, r\sin\theta, z):=q(\theta,r,z)$. For a fixed angle $\bar{\theta}$, $q$ is
a polynomial of degree at most $n$ in $r,z$, while it is a trigonometric polynomial in $\theta$ of degree at most $n$ for fixed values of $(r,z)$. Since on the generic rectangle $(r,z)$ we have considered the set of Padua points of degree $n$ which
is a WAM, say $A_1$, hence, we can write
$$|q(\bar{\theta},r,z)| \le c_1(n) \|q(\bar{\theta},\cdot,\cdot)\|_{A_1}, \;\; $$
where $c_1(n)$ does not depend on $\bar{\theta}$. Let $|q(\bar{\theta},r^*,z^*)|$ be the maximum.
Considering now the equispaced angles, 
$\theta_k=2k\pi/(2n+2),\; 0\le k\le 2n+1$
i.e. the $2n+2$ equispaced points in $[0,2\pi[$, say $A_2$, then
$$|q(\bar{\theta},r^*,z^*)|\le c_2(n) \|q(\cdot,r^*,z^*)\|_{A_2}\,.$$

Passing to the maximum also on the left side, we have
$$ \max_{(x,y,z) \in K} |p(x,y,z)|=\|p\|_K \le c_1(n)c_2(n) \|q\|_{A_n}$$
that is $\|p\|_K \le C(n)\|p\|_{A_n}$, where $C(n)$ is indeed $C(A_n)={\cal O}(\log^3 n)$, showing that this discretization is a WAM for the cylinder $K$.
\vskip 0.2in

\noindent
\begin{figure}[hbt]
  \begin{center}
    \begin{tabular}{c}
    \includegraphics[scale=0.5]{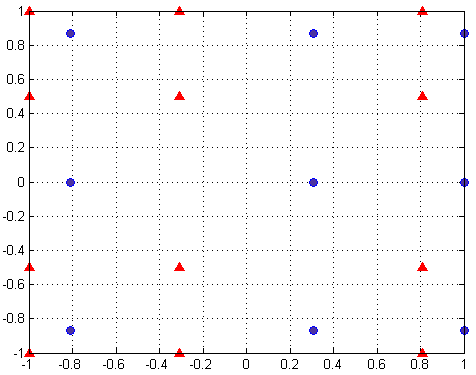}
      \includegraphics[scale=0.5]{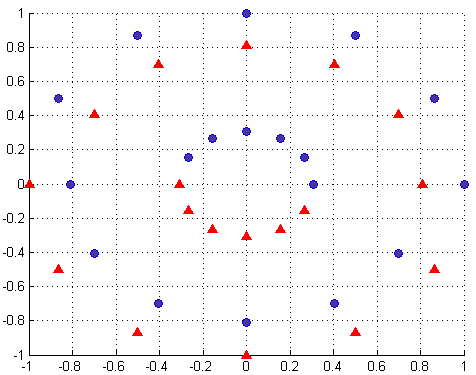}
      \end{tabular}
    \caption{
Left: Padua points for $n=5$ on the $(r,z)$-plane. 
We used different markers for the two Chebsyhev meshes of Padua points. 
Right: Projections on the disk of WAM2 points  
(obtained by $n+1$ rotations of Padua points around the $z$ axis). Notice that the points
indicated with the small triangles are on different levels. See also Figure \cite{Fig1}.} \label{Fig0}
    \end{center}
\end{figure}

\begin{figure}[hbt]
  \begin{center}
    \begin{tabular}{c}
      \includegraphics[scale=0.5]{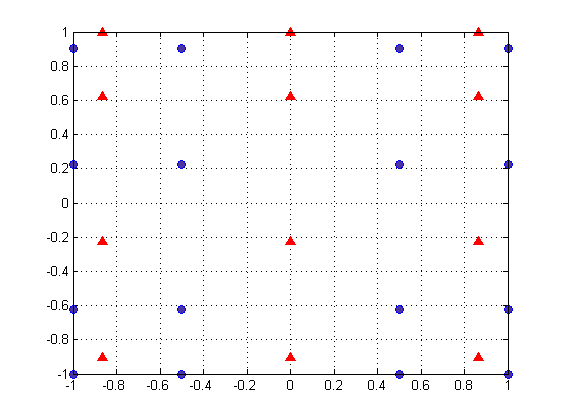}
      \includegraphics[scale=0.5]{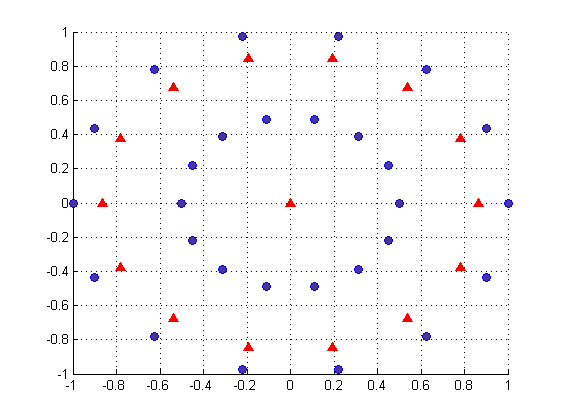}
      \end{tabular}
    \caption{
Left: Padua points for $n=6$ on the $(r,z)$-plane. 
We used different markers for the two Chebsyhev meshes of Padua points. 
Right: Projections on the disk of WAM2 points  
(obtained by $n+1$ rotations of Padua points around the $z$ axis). Notice that the points
indicated with the small triangles are on a different levels.} \label{Fig01}
    \end{center}
\end{figure}

\begin{figure}
  \begin{center}
    \begin{tabular}{c}
      \includegraphics[scale=0.8]{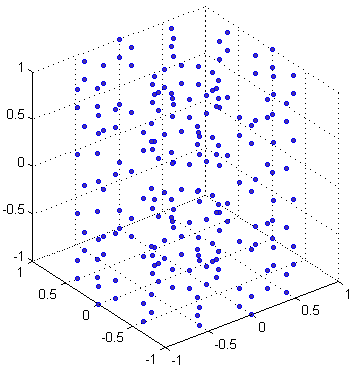}\\
      \includegraphics[scale=0.8]{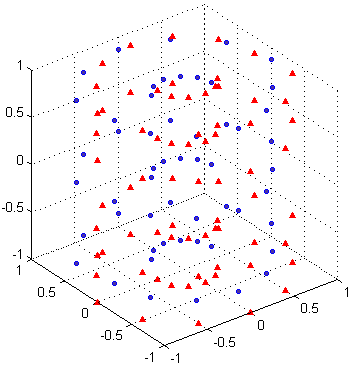}
      \includegraphics[scale=0.8]{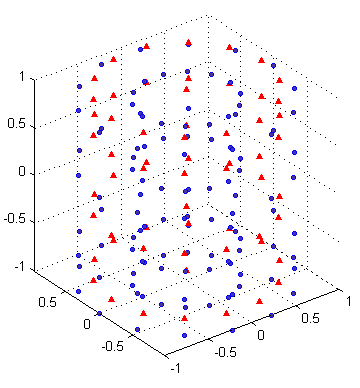}
      \end{tabular}
    \caption{
Above: the {\em first WAM} for $n=5$ having $216$ points.
Below left: the {\em second WAM} for $n=5$ having $126$ points. 
Below right: the {\em second WAM} for $n=6$ having $172$ points.} \label{Fig1}
    \end{center}
\end{figure}

\section{Computation of AFP and DLP}
As discussed in \cite{BDMSV09}, the computation the AFP and DLP, can
be done by a few basic linear algebra operations, corresponding to
the LU factorization with row pivoting of the Vandermonde matrix for
the DLP, and to the QR factorization with column pivoting of the
transposed Vandermonde matrix for the AFP (cf. \cite{SV09-1}). For
the sake of completeness, we recall these two Matlab-like scripts
used in \cite{SV09-1,BDMSV09} for computing the AFP and DLP,
respectively. \vskip0.1cm \noindent
{\bf algorithm AFP (Approximate Fekete Points):}\\
$\bullet$ $W=(V(\bs{a},\bs{p}))^t$; $\bs{b}=(1,\dots,1)^t\in
\mathbb{C}^N$; $\bs{w}=W\backslash \bs{b}\,$;
$ind=\mbox{\texttt{find}} (\bs{w}\neq \bs{0})$;
$\bs{\xi}=\bs{a}(ind)$ \vskip0.1cm \noindent \vskip0.1cm \noindent
{\bf algorithm DLP (Discrete Leja Points):}\\
$\bullet$ $V=V(\bs{a},\bs{p})$;
$[L,U,\bs{\sigma}]=\mbox{\texttt{LU}}(V,\mbox{``\texttt{vector}''})$;
$ind=\bs{\sigma}(1,\dots,N)$; $\bs{\xi}=\bs{a}(ind)$ \vskip0.1cm

\begin{figure}
  \begin{center}
    \begin{tabular}{c}
      \includegraphics[scale=0.7]{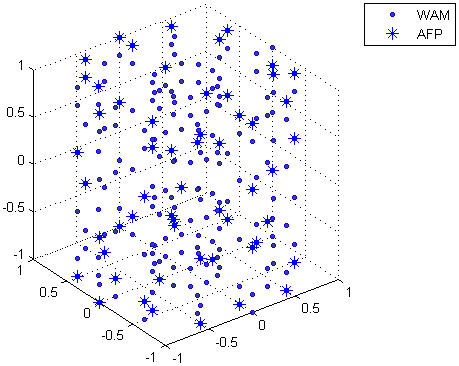}
      \includegraphics[scale=0.7]{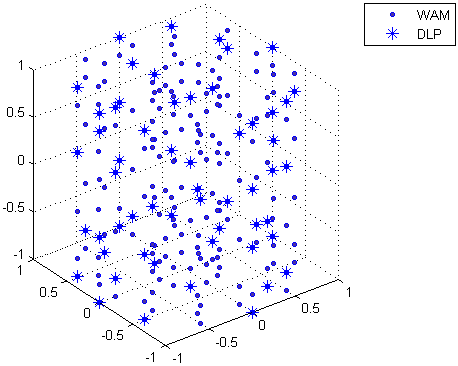}
      \end{tabular}
    \caption{
WAM1 of the cylinder for $n=5$ and the corresponding extracted points.
Left: $56$ Approximate Fekete Points. Right: $56$ Discrete Leja Points.} \label{Fig2}
    \end{center}
\end{figure}

\begin{figure}
  \begin{center}
    \begin{tabular}{c}
      \includegraphics[scale=0.7]{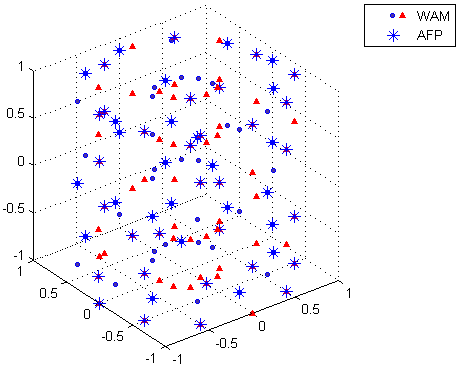}
      \includegraphics[scale=0.7]{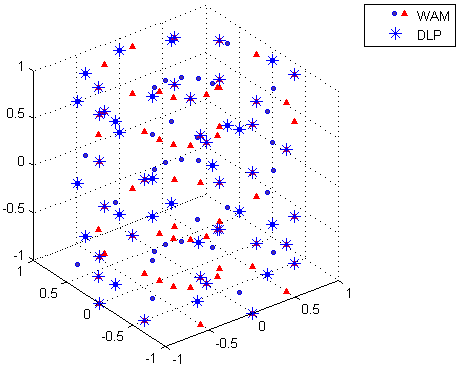}
      \end{tabular}
    \caption{
WAM2 of the cylinder for $n=5$ and the corresponding extracted
points. Left: $56$ Approximate Fekete Points. Right: $56$ Discrete
Leja Points.} \label{Fig3}
    \end{center}
\end{figure}

In Figures \ref{Fig2}--\ref{Fig3}, we show the AFP and DLP extracted
from the WAM1 and WAM2 for $n=5$. \vskip 0.1in \noindent In the
above scripts, $V(\bs{a},\bs{p})$ indicates the Vandermonde matrix
at the WAM $\bs{a}$ using the polynomial basis $\bs{p}$, that is the
matrix whose elements are $p_j(a_i), \;1\le j\le N,\; 1\le i\le {\tt
card}(A_n)$. The extracted AFP and DLP are then stored in the vector
$\bs{\xi}$. \vskip 0.1in \noindent {\bf Remark 1}. In both
algorithms, the selected points (as opposed to the continuum Fekete
points) depend on the choice of the polynomial basis. But in the
second algorithm, which is based on the notion of determinant (as
described in \cite[\S6.1]{BDMSV09}), the selected points also depend
on the ordering of the basis. In the univariate case with the
standard monomial basis, it is not difficult to recognize that the
selected points are indeed the Leja points extracted from the mesh
(cf. \cite{BC09,ST97} and references therein). \vskip 0.1in
\noindent {\bf Remark 2}. When the conditioning of the Vandermonde
matrices is too high, and this happens when the polynomial basis is
ill-conditioned, the algorithms can still be used provided that a
preliminary iterated orthogonalization, that is a change to a
discrete orthogonal basis, is performed (cf.
\cite{BCLSV09,BDMSV09,SV09-1}). This procedure however only
mitigates the effect of a bad choice of the polynomial basis.
Consequently, whenever is possible, is desirable to use a
well-conditioned polynomial basis.

\vskip 0.2in \noindent
A suitable basis for the forementioned rectangular cylinder $K$,
is the set of polynomials introduced by J. Wade in \cite{W10}:
\begin{equation} \label{baseWade}
C_{j,k,i}(x,y,z):=U_k(\theta_{j,k};x,y)\tilde{T}_{i-k}(z),\; \; i=0,\ldots,n,\,\, k=0,\ldots,i,\,\, j=0,\ldots,k
\end{equation}
where
\begin{itemize}
\item $\theta_{j,k}=\frac{j\pi}{k+1}$;
\item $U_k(\theta_{j,k};x,y)=U_k(x\cos(\theta_{j,k})+y\sin(\theta_{j,k}))$ is the Chebyshev polynomial of the second kind which is
an orthonormal basis for the disk w.r.t. the measure
$\omega(x,y)=\frac{1}{\pi}$;
\item $\tilde{T}_{j}(z)$ is the $j$-th orthonormal Chebyshev polynomial of the first kind,
i.e. $\tilde{T}_{j}(z)=\sqrt{2}T_{j}(z)$, w.r.t. the measure
$\omega(z)=(1-z^2)^{-\frac{1}{2}}$.
\end{itemize}

As discussed in \cite{W10}, the basis $C_{j,k,i}(x,y,z)$ is
orthonormal for the space of orthogonal polynomials on $K$ w.r.t.
the weight function $\frac{1}{\pi}{(1-z^2)}^{-\frac{1}{2}}$. This
basis plays also an important role in the construction of the
discretized Fourier orthogonal expansions on the disk $D$ and the
unitary rectangular cylinder. Moreover, this turns out to be
well-conditioned. Indeed, as outlined in \cite{W10}, if we consider
the Radon projection of a function $f$ on $\mathbb{R}^d$
$$ R_\theta(f;t)=\int_{<x,\theta>=t} f(x)dx,\;\; t \in \mathbb{R},\; \theta \in \mathbb{S}^{d-1}$$
with $\mathbb{S}^{d-1}$ the unit sphere in $\mathbb{R}^d$.
The discretized Fourier expansion $F_{2m}$ of $R_\theta(f;t)$ for $m\ge 0$, $(x,y)$ belonging to the unit disk and $z \in [-1,1]$
$$
F_{2m}(f)(x,y,z)=\sum_{i=0}^{2m}\sum_{j=1}^{2m}\sum_{k=0}^{2m-1}
R_{\phi_\nu}(f(\cdot,\cdot,z_{k,2m}; \cos(\theta_{j,2m}))) {\cal
T}_{i,j,k}(x,y,z)
$$
where $\phi_\nu={2\nu \pi \over 2m+1}$,
$$ {\cal T}_{i,j,k}(x,y,z)={1 \over 4m^2(2m+1)}\sum_{n=0}^{2m}\sum_{l=0}^{n} (l+1)\sin((l+1)\theta_{j,2m}) U_l(\cos(\sigma_i(x,y)))T_{n-l}(z_{k,2m})T_{n-l}(z)
$$
and $\sigma_i(x,y)=\arccos(x \cos \phi_i + y \sin \phi_i)$, has sup norm that
grows as $\|F_{2m}\|_\infty \approx m (\log(m+1)^2)$, i.e. nearly the optimal growth.

\section{Approximation and cubature on the cylinder}

\subsection{Interpolation and least-square approximation}
The interpolation polynomial $q_n({\bf x})$ of degree $n$ of a real
continuous function $f$ defined on the compact $K \subset
\mathbb{R}^3$, can be written in Lagrange form as
\begin{equation}
q_n({\bf x})=\sum_{j=1}^N f({\bf a}_j) l_j({\bf x}), \; {\bf x} \in
K
\end{equation}
where $N={\tt dim}(\mathbb{P}_n^3)$, ${\bf a}_j$ are the AFP or the
DLP extracted from the WAMs and $l_j$ indicates the $j$th elementary
Lagrange polynomial of degree $n$. Let ${\bf l}=(l_1({\bf
x}),\ldots,l_N({\bf x}))$ be the (row) vector of all the elementary
Lagrange polynomials at a point ${\bf x}$, ${\bf p}$ the vector of
the basis (\ref{baseWade}) and ${\bf a}=({\bf a}_1,\ldots, {\bf
a}_N)$, then we can compute ${\bf l}$ by solving the linear system
$$ {\bf l}^t=W \,{\bf p}^t,\;\; W=(V({\bf a}, {\bf p})^{-1})^t\,.$$
The interpolation operator ${\cal L}_n: {\cal C}(K) \rightarrow
\mathbb{P}^3_n$, with ${\cal C}(K)$ equipped with the sup norm, that
maps every $f \in {\cal C}(K)$ into the corresponding polynomial $q
= \sum_{i=1}^N l_i(\cdot)f({\bf x}_i) \in \mathbb{P}^3_n$ in
Lagrange form, is a projection having norm
\begin{equation}
\|{\cal L}_n\|_\infty=\max_{{\bf x} \in K} \sum_{j=1}^N|l_j({\bf x})|:= \Lambda_n
\end{equation}
where $\Lambda_n$ is the well-known {\em Lebesgue constant}.
When the interpolation points are the true Fekete points,
the Lebesgue constant satisfies the upper bound
$$\Lambda_n =\max_{{\bf x} \in K} \|W\,{\bf p}^t\|_1\le N\,$$
since $\|l_j\|\le 1$.

Thanks to property {\bf P10} of WAMs we can say more. Indeed, when
the Fekete points are extracted from a WAM, $\|l_j\|_K \le C(A_n)
\|l_j\|_{A_n} \; \forall j$ (cf. \cite[\S4.4]{CL08}), 
from which the following upper bound
holds
$$\Lambda_n \le N\, C(A_n)\,,$$
with $C(A_n)$, the same constant in definition of a WAM, 
which depends on $A_n$. In the
numerical experiments that we will present in the next section, we
will observe that the above upper bound is a quite pessimistic
overestimate. \vskip 0.1in \noindent Another natural application of
such a construction, is the {\em least squares approximation} of a
function $f \in {\cal C}(K)$. Given a WAM $A_n=\{{\bf a}_1,
\ldots,{\bf a}_M\}, \, M \ge N$, with $N=\mbox{dim}(\mathbb{P}_n^3)$
and ${\bf p}=\{p_1,\ldots,p_N\}$ a basis for $\mathbb{P}_n^3$, let
us consider the orthonormal basis ${\bf q}$ w.r.t. the discrete
inner product $\left<f,g\right>=\sum_{i=1}^M f({\bf
a}_i)\overline{g({\bf a}_i)}$ which can be obtained from the basis
${\bf p}$ by multiplying by a certain transformation matrix $P$,
i.e. {\bf q}=P{\bf p}.

The least squares operator of $f$ at the points of a WAM $A_n$ can
then be written as
\[
{L}_{A_n}(f)({\bf x})= \sum_{j=1}^N\left( \sum_{i=1}^M f({\bf a}_i) \overline{q_j({\bf a}_i)}\right)
q_j({\bf x})=\sum_{i=1}^M f({\bf a}_i) g_i({\bf x})
\]
where $g_i({\bf x})=\sum_{j=1}^Nq_j({\bf x}) \overline{q_j({\bf
a}_i)}$. Letting ${\bf g}=(g_1,\ldots,g_M)^t$, it follows that ${\bf
g}=Q P^t{\bf p}$, where the matrix $Q$ is a numerically orthogonal
(unitary) matrix, i.e. $\bar{Q^t}Q=I$, and $Q=V({\bf a}, {\bf
q})=V({\bf a}, {\bf p})P$. Notice that, the transformation matrix
$P$ and the matrix $Q$ are computed once and for all for a fixed
mesh.

\noindent The norm of the operator, that is its Lebesgue constant,
is then
\[
 \|{L}_{A_n}\|=\max_{{\bf x}\in K} \sum_{i=1}^M |g_i({\bf x})|=\max_{{\bf x}\in K} \|QP^t{\bf p}({\bf x})\|_1\,.
 \]
In \cite{CL08}, it is observed that three-dimensional WAMs so defined can be used as discretization
of compact sets $K \subset \mathbb{R}^3$ for the computation of good points for least-square approximation
by polynomials. Indeed, using property {\bf P9}, in \cite[Th. 2]{CL08} the authors proved the following
error estimates for least-squares approximation on WAMs of a function $f \in {\cal C}(K)$
\begin{equation} \label{lserror}
\|f -L_{A_n}(f)\|_K \le \left(1+C(A_n)(1+\sqrt{{\tt
card}(A_n)})\right)\min\{\|f -p\|_K\,:\,p \in \mathbb{P}_n^3\}
\end{equation}
where again $C(A_n)$ depends on the WAM.

This estimates says that if we could control the factor
$C(A_n)(1+\sqrt{{\tt card}(A_n)})$, and this is possible when we use
WAMs, then the approximation $L_{A_n}(f)\in \mathbb{P}_n^3 $ is
nearly optimal.

\subsection{Cubature}
For a given function $f: K\subset \mathbb{R}^3 \rightarrow
\mathbb{R}$ we want to compute
\[
I(f)=\int_K f({\bf x}) d{\bf x}
\]
where $d{\bf x}$ is the usual Lebesgue measure of the compact set
$K$. An interpolating cubature formula $C_N(f)$ that approximates
$I(f)$ can be expressed as
$$C_N(f)=\sum_{i=1}^N w_i f({\bf x}_i)\,$$
where, in our case, the nodes ${\bf x}_i$ are the AFP or the DLP for
$K$. Once we know the {\em cubature weights} $w_i$, the $C_N(f)$
gives an approximation of $I(f)$. The cubature weights can be
determined by solving the {\em moment system}, that is
\[
\sum_{j=1}^N w_j p_i({\bf x}_j)=\int_K p_i({\bf x})d{\bf x}, \;\; i=1,\ldots,N
\]
where $p_i$ is the $i$-th element of the polynomial basis {\bf p}.
Hence, if $V=(p_i(x_j))$ and $b_i=\int_K p_i(x)dx$, the nodes $x_i$
and the weights $w_i$ are provided by the AFP or DLP algorithm.

\section{Numerical results}
In this section we present the numerical experiments that we made for showing the quality of AFP and DLP on interpolation,
approximation by least-squares
and cubature of $6$ test functions defined on the rectangular cylinder $K=D\times [-1,1]$.
The results are obtained by using both the AFP and the DLP extracted from both WAM1 and WAM2.

In Table \ref{Tab1} we collect, for degrees $n=5,10,15,20,25,30$, the values of the Lebesgue constant $\Lambda_n$,
those of the condition number (in the sup norm) $\kappa_{\infty,1}$ for the
associated Vandermonde matrix using the Wade basis for the WAM1 for the AFP.
Due to hardware restrictions, the Lebesgue constant has been evaluated on a control mesh $M_n=A_{m}$ with $m=4n$ for $n\le 20$,
$m=2n$ for $n > 20$.

In Table \ref{Tab2} we present $\Lambda_n, \, \kappa_{\infty,2}$ for
the WAM2 on the AFP. In this case the control mesh $M_n=A_m$ with
$m=4n$ for $n\le 20$, $m=3n$ for $n \le 25$ and $m=2n$ for $n> 25$.

In Tables \ref{Tab3} and \ref{Tab4} we present the corresponding values for the DLP.

In Table \ref{Tab5} we display the the norm of the least-square operator
$\|L_{A_n}\|$ on the WAM1 and WAM2, respectively.

Note that the norm of the least-square operator $\|L_{A_n}\|$ (norm
computed at the points of the WAM $A_n$) has been computed after two
steps of orthonormalization of the polynomial basis.

The numerical results show that both $\Lambda_n$ and the condition
number $\kappa_{\infty,*}$ are smaller for WAM2 while, by oppositeon
the contrary, the sup norm of the least-square operator turns out to
be bigger for WAM2 than WAM1. Actually, for WAM2, $\|L_{A_n}\|$ has
a growth factor close to $2$.

\begin{table}[hbt]
\begin{center}
\begin{tabular}{cccccccc}
\hline
n & 5 & 10 & 15 & 20 & 25 & 30 \\
\hline
\footnotesize{$\Lambda_n$}&17 &83 &208 &384 &849 & 988\\
\footnotesize{$\kappa_{\infty,1}$}&18.2 &177 &384 &746 &1410 & 2650\\
\hline
\end{tabular}
\end{center}
\caption{The Lebesgue constant and the condition number of the Vandermonde matrix constructed
using the Wade basis on the AFP for the WAM1.} \label{Tab1}
\end{table}

\begin{table}[hbt]
\begin{center}
\begin{tabular}{cccccccc}
\hline
n & 5 & 10 & 15 & 20 & 25 & 30 \\
\hline
\footnotesize{$\Lambda_n$} &19 &76 &213 &427 &879 &1034 \\
\footnotesize{$\kappa_{\infty,2}$}&19.4 &115 &440 &705 &1540 &2380\\
\hline
\end{tabular}
\end{center}
\caption{The Lebesgue constant and the condition number of the
Vandermonde matrix constructed using the Wade basis on the AFP for the WAM2. } \label{Tab2}
\end{table}

\begin{table}[hbt]
\begin{center}
\begin{tabular}{cccccccc}
\hline
n & 5 & 10 & 15 & 20 & 25 & 30 \\
\hline
\footnotesize{$\Lambda_n$}&30 &115 &350 &617 &1388 &2597\\
\footnotesize{$\kappa_{\infty,1}$}&35.2 &247 &772 &1190 &3090 &5320\\
\hline
\end{tabular}
\end{center}
\caption{The Lebesgue constant and the condition number of the Vandermonde matrix constructed
using the Wade basis on the DLP for the WAM1.} \label{Tab3}
\end{table}

\begin{table}[hbt]
\begin{center}
\begin{tabular}{cccccccc}
\hline
n & 5 & 10 & 15 & 20 & 25 & 30 \\
\hline
\footnotesize{$\Lambda_n$} &30 &129 &349 &648 &1520 &2143 \\
\footnotesize{$\kappa_{\infty,2}$}&29.9 &176 &638 &782 &2310 &3940\\
\hline
\end{tabular}
\end{center}
\caption{The Lebesgue constant and the condition number of the
Vandermonde matrix constructed using the Wade basis on the DLP for the WAM2.} \label{Tab4}
\end{table}

\begin{table}[hbt]
\begin{center}
\begin{tabular}{cccccccc}
\hline
n & 5 & 10 & 15 & 20 & 25 & 30 \\
\hline
\footnotesize{$||{L}_{A_{n,1}}||$} &4.8& 10.2& 10.7& 21.1& 15.8& 22.6\\
\footnotesize{$||{L}_{A_{n,2}}||$} &7.2& 15.3& 32.8& 43.4& 85.9& 96.6\\
\hline
\end{tabular}
\end{center}
\caption{The sup-norm of the least-squares operator on WAM1 and WAM2, respectively. } \label{Tab5}
\end{table}
\newpage
The interpolation and the cubature relative errors on the AFP and the DLP have been computed for
the following six test functions:
\begin{eqnarray*}
f_1(x,y,z) &=&      0.75{\tt e}^{-\frac{(9x-2)^2+(9y-2)^2+(9z-2)^2}{4}}+0.75{\tt e}^{-\frac{(9x+1)^2}{49}-\frac{9y+1}{10}-\frac{9z+1}{10}}\\
&+& 0.5{\tt e}^{-\frac{(9x-7)^2+(9y-3)^2+(9z-5)^2}{4}}-0.2{\tt e}^{-(9x-4)^2-(9y-7)^2-(9z-5)^2};\\
f_2(x,y,z) &=& \sqrt{(x-0.4)^2+(y-0.4)^2+(z-0.4)^2};\\
f_3(x,y,z) &=& \cos(4(x+y+z));\\
f_4(x,y,z) &=& \frac{1}{1+16(x^2+y^2+z^2)};\\
f_5(x,y,z) &=& \sqrt{(x^2+y^2+z^2)^3};\\
f_6(x,y,z) &=& \cos(x^2+y^2+z^2)\,.
\end{eqnarray*}
The function $f_1$ is the three-dimensional equivalent of the well-known {\em Franke test function}.
The function $f_2$ has a singular point into the cylinder $K=D\times [-1,1]$.
The functions $f_3$ and $f_6$ are infinitely differentiable. The function $f_4$ is the {\em Runge function}. The function $f_5$ is
a ${\cal C}^2$ function with third derivatives singular at the origin.
\vskip 0.15in

All numerical experiments have been done on a cluster HP with 14
nodes. We used one of the nodes equipped with 2 processors quad core
with 64Gb of RAM.

\vskip 0.15in

\begin{figure}[hbt]
  \begin{center}
    \begin{tabular}{cc}
      \resizebox{75mm}{!}{\includegraphics{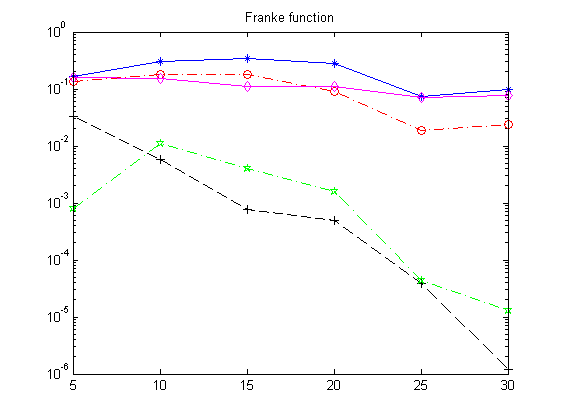}} &
      \resizebox{75mm}{!}{\includegraphics{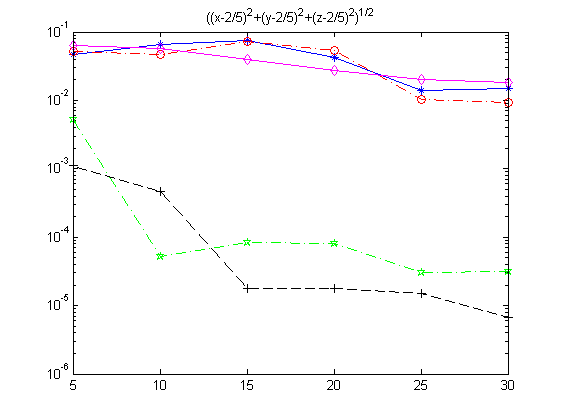}} \\
      \resizebox{75mm}{!}{\includegraphics{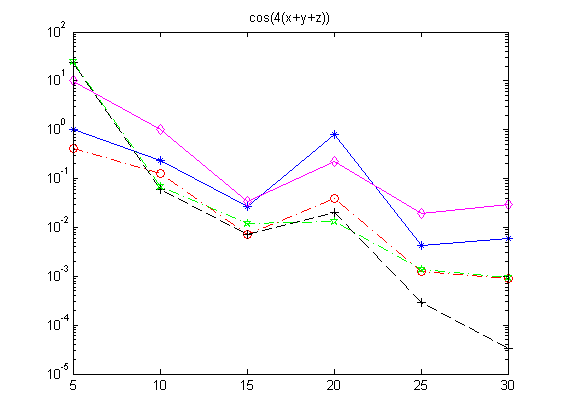}} &
      \resizebox{75mm}{!}{\includegraphics{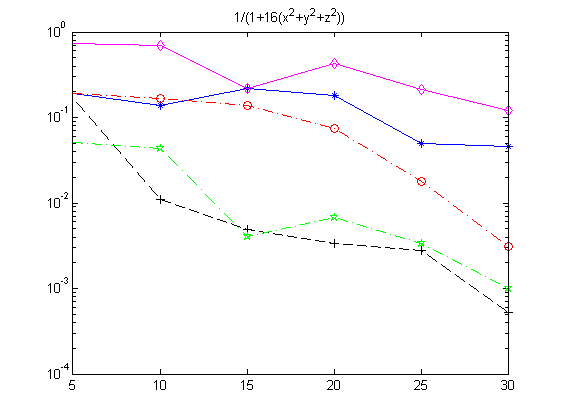}} \\
      \resizebox{75mm}{!}{\includegraphics{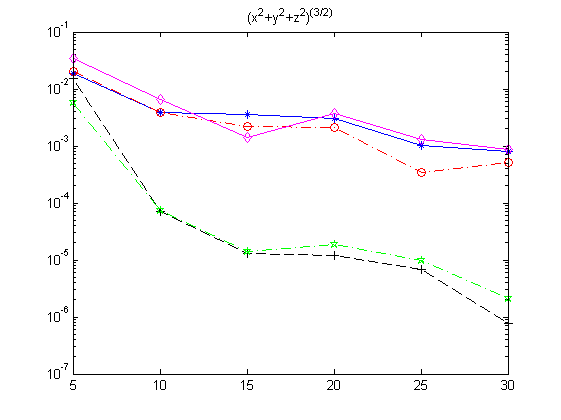}} &
      \resizebox{75mm}{!}{\includegraphics{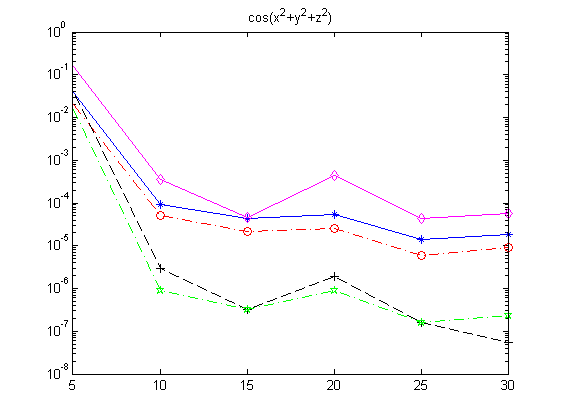}} \\
    \end{tabular}
    \caption{Interpolation error on AFP (\textcolor{blue}{blue $-*$}), cubature errors on AFP (\textcolor{black}{black $-+$})
    Interpolation error on DLP (\textcolor{magenta}{magenta $-\diamond$}), cubature errors on DLP (\textcolor{green}{green $--\star$})
    and least-squares errors (\textcolor{red}{red $-\cdot$o}).
    The points are extracted from the WAM1. Wade basis. In abscissa the polynomial degree.}\label{cili1err1}
    \end{center}
\end{figure}

\begin{figure}[hbt]
  \begin{center}
    \begin{tabular}{cc}
      \resizebox{75mm}{!}{\includegraphics{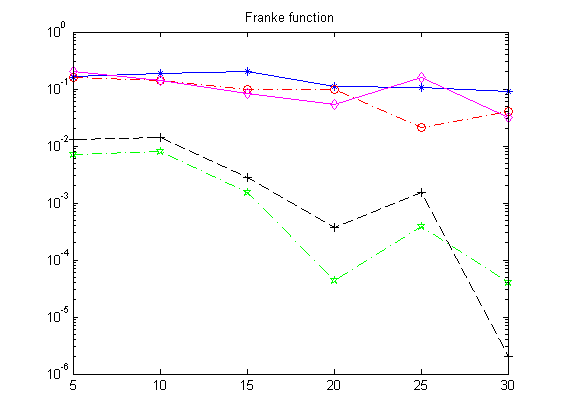}} &
      \resizebox{75mm}{!}{\includegraphics{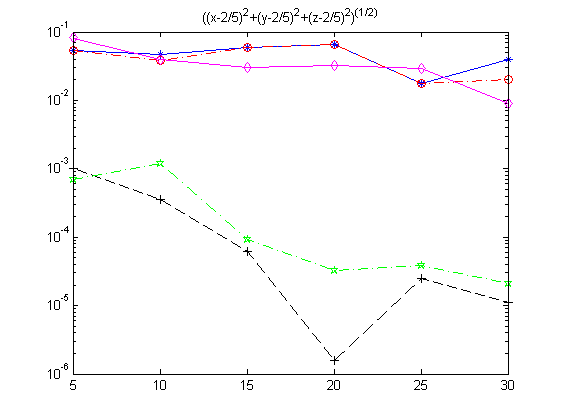}} \\
      \resizebox{75mm}{!}{\includegraphics{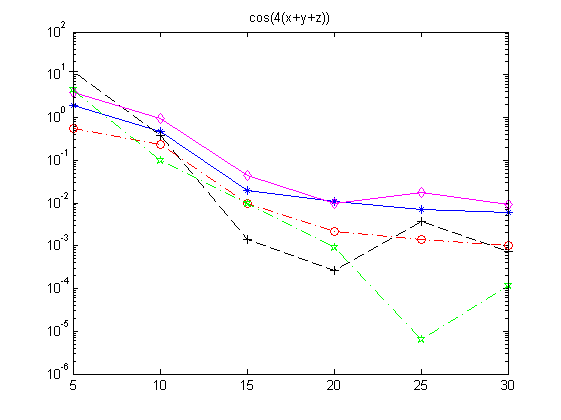}} &
      \resizebox{75mm}{!}{\includegraphics{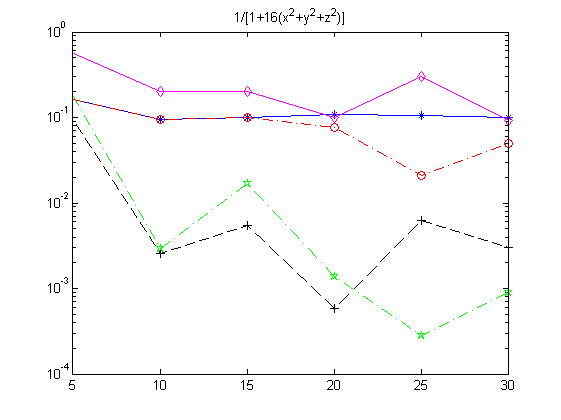}} \\
      \resizebox{75mm}{!}{\includegraphics{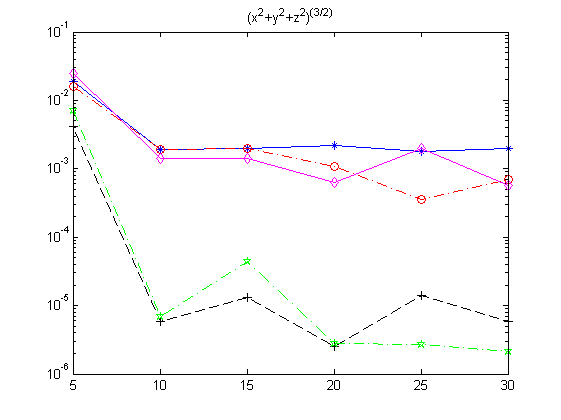}} &
      \resizebox{75mm}{!}{\includegraphics{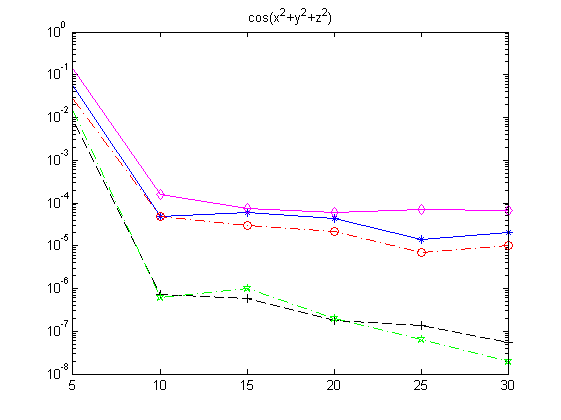}} \\
    \end{tabular}
    \caption{Interpolation error on AFP (\textcolor{blue}{blue $-*$}), cubature errors on AFP (\textcolor{black}{black $-+$})
    Interpolation error on DLP (\textcolor{magenta}{magenta $-\diamond$}), cubature errors on DLP (\textcolor{green}{green $--\star$})
    and least-squares errors (\textcolor{red}{red $-\cdot$o}).
    The points are extracted from the WAM2. Wade basis. In abscissa the polynomial degree.}\label{cili1err2}
    \end{center}
\end{figure}

In Figures \ref{cili1err1} and \ref{cili1err2} we display the
interpolation, cubature and least-square relative errors on the AFP
and DLP, up to degree $n=30$, for the WAM1 and WAM2, respectively.
The results shows that the AFP give, in general, smaller errors.
Only the cubature errors on WAM2 are smaller for DLP than AFP. One
reason is related to the values of the Lebesgue constants and the
conditioning of the Vandermonde matrices that are smaller for AFP
than DLP.

Since WAM2 has a lower cardinality and that the results are more or
less the same, such a mesh is more convenient from the point of view
of efficiency and approximation order. As true values of the
functions $f_i, \; i=1,\ldots,6$, we considered the value of $f_i$
on the control meshes used for computing the Lebesgue constants. As
exact values of the integrals, we considered the values computed by
the Matlab built-in function {\tt triplequad} with the chosen
tolerance depending on the smoothness of the function. For the
smoothest functions $f_3$ and $f_6$ we used the tolerance $1.e-12$
while for the others $1.e-10$. This choice allowed to avoid stalling
phenomena that we encountered in computing the integrals of $f_3$
and $f_6$.

We point out that in the case of the function $f_2$, where there
exist a singularity in $(0.4,0.4,0.4)$, {\tt triplequad} uses a
domain decomposition approach, implemented in the method {\tt
quadgk} (Gauss-Kronrod cubature rules) avoiding the singularity.
Actually, due to the geometry of the cylinder, we could compute the
exact values for the integrals by using separation of variables,
that is instead of a call to the Matlab built-in function {\tt
triplequad} we used the product of the built-in functions {\tt
dblquad} and {\tt quadl}. This allowed a considerably reduction of
the computational time, as displayed in Table \ref{Tab6} for degrees
$n \le 20$.

Concerning function $f_4$, the least-square errors seem not those that
one can expected. In \cite[Fig. 3.2]{DMVX09} the authors already computed the
relative hyperinterpolation errors for the Runge function w.r.t. the number
of function evaluations. From that figure, correspondingly to polynomial
degree $n=30$, that requires $(n+1)(n+2)(n+3)/6=5456$ function evaluations,
the hyperinterpolation error is about $10^{-1}$. Hence, what we see 
in Figure \ref{cili1err2} is consistent with those results and, as expected, 
formula (\ref{lserror}) is an overestimate of the least-square error.
\begin{table}[hbt]
\begin{center}
\begin{tabular}{c|l|l}
\hline
n & {\tt triplequad} & {\tt dblquad}$\times${\tt quadl} \\
\hline
5 & 1 min. 9 sec. & 4.7 sec.\\
10& 16 min. 3 sec. & 1 min. 4 sec.\\
15& 1 h 15 min. 10 sec. & 5 min. 33 sec. \\
20& 3 h 47 min. 12 sec. & 16 min. 30 sec.\\
\hline
\end{tabular}
\end{center}
\caption{Computational time: {\tt triplequad} vs {\tt
doublequad}$\times${\tt quadl}} \label{Tab6}
\end{table}

\vskip 0.2in \noindent {\bf Acknowledgments}. This work has been
done with the support of the {\it 60\% funds, year 2010} of the
University of Padua.


\bibliographystyle{plain}

\begin{thebibliography}{99}

\bibitem{BC09} L.~Bialas-Cie{\'z} and J.-P.~Calvi, {\em Pseudo Leja
sequences}, Ann. Mat. Pura Appl., published online November 16,
2010.

\bibitem{BCDMVX06} L.~Bos, M. Caliari, S. De Marchi, M. Vianello and Y. Xu,
{\em Bivariate Lagrange interpolation at the Padua points: the generating curve approach},
J. Approx. Theory 143 (2006), 15--25.

\bibitem{BCLSV09} L.~Bos, J.-P.~Calvi, N.~Levenberg, A.~Sommariva
and M.~Vianello, {\em Geometric weakly admissible meshes, discrete
least squares approximation and approximate Fekete points}, to
appear in Math. Comp. (2010) (preprint online at:
http://www.math.unipd.it/$\sim$marcov/CAApubl.html).

\bibitem{BDMSV09} L.~Bos, S.~De Marchi, A.~Sommariva and M.~Vianello,
{\em Computing multivariate Fekete and Leja points by numerical
linear algebra}, SIAM J. Numer. Anal. 48 (2010), 1984--1999.

\bibitem{BDMSV09_1} L.~Bos, S.~De Marchi, A.~Sommariva and M.~Vianello,
{\em Weakly Admissible Meshes and Discrete Extremal Sets}, Numer.
Math. Theory Methods Appl. 4 (2011), 1--12.

\bibitem{BSV09} L.~Bos, A.~Sommariva and M.~Vianello,
{\em Least-squares polynomial approximation on weakly admissible
meshes: disk and triangle}, J. Comput. Appl. Math. 235 (2010), 660--668.

\bibitem{DMVX09} S. De Marchi, M. Vianello and Y. Xu,
{\em New cubature formulae and hyperinterpolation in three variables},
BIT, Vol. 49(1) 2009, 55-73.

\bibitem{CL08} J. P.~Calvi and N.~Levenberg, {\em Uniform
approximation by discrete least squares polynomials}, J. Approx.
Theory 152 (2008), 82--100.

\bibitem{CMI09} A.~Civril and M.~Magdon-Ismail, {\em On selecting
a maximum volume sub-matrix of a matrix and related problems},
Theoretical Computer Science 410 (2009), 4801--4811.

\bibitem{KL92} M. Klimek, Pluripotential Theory, Oxford U. Press,
1992.

\bibitem{K10} A. Kro\'o, {\em On optimal polynomial meshes},
J. Approx. Theory, to appear.

\bibitem{Tesi} M. Marchioro, {\em Punti di Approssimati di Fekete e Discreti di Leja del parallelepipedo,
del cilindro e del prisma retto}, Master Thesis, University of Padua
(2010) (in Italian).

\bibitem{PV10} F.~Piazzon and M. Vianello, {\em Analytic transformations of admissible meshes},
East J. Approx. 16 (2010), 313--322.

\bibitem{ST97} E.B.~Saff and V.~Totik, Logarithmic potentials with
external fields, Springer, 1997.

\bibitem{SV09-1} A.~Sommariva and M.~Vianello, {\em Computing
approximate Fekete points by QR factorizations of Vandermonde
matrices}, Comput. Math. Appl. 57 (2009), 1324--1336.

\bibitem{SV09-2} A.~Sommariva and M.~Vianello, {\em Gauss-Green cubature and
moment computation over arbitrary geometries}, J. Comput. Appl.
Math. 231 (2009), 886--896.

\bibitem{SV09-3} A.~Sommariva and M.~Vianello, {\em Approximate
Fekete points for weighted polynomial interpolation}, Electron.
Trans. Numer. Anal. 37 (2010), 1--22.

\bibitem{W10} J. Wade, {\em A discretized Fourier orthogonal expansion in orthogonal polynomials on a
cylinder}, J. Approx. Theory 162 (2010), 1545--1576.

\end{thebibliography}

\end{document}